\documentclass{article}

\usepackage{enumerate} 
\usepackage{amsmath}  %
\usepackage{amsfonts} 
\usepackage{amssymb}  %
\usepackage{amsthm}   %

\newcommand{\RR}[0]{\mathbb{R}}  
\newcommand{\CC}[0]{\mathbb{C}}  

\newcommand{\pd}[0]{\partial} 

\newcommand{\frk}[1]{\mathfrak{#1}} 
\newcommand{\cc}[1]{\overline{#1}} 
\newcommand{\nin}{\notin} 
\newcommand{\mc}[1]{\mathcal{#1}} 
\newcommand{\bksl}{\backslash}

\theoremstyle{plain}

\newtheoremstyle{named}{}{}{\itshape}{}{\bfseries}{.}{.5em}{#2. \thmnote{#3}}
\theoremstyle{named}

\newcommand{\TSpec}{\sigma_{\mathrm{Ta}}}

\title{The Taylor joint spectrum and restriction to hyperinvariant subspaces}
\author{Edward J. Timko\thanks{Partially supported by a Pacific Institute of Mathematical Sciences Fellowship}}

\begin{document}
\maketitle

\begin{abstract}
	It is well known that for a single bounded operator $A_0$ on a Hilbert $\frk{H}$, if $\frk{M}\subset \frk{H}$ is hyperinvariant for $A_0$, then the spectrum of $A_0|_{\frk{M}}$ is contained in the spectrum of $A_0$.
	In this note, we modify an example of Taylor to prove the following.
	There exist a quadruple $A=(A_1,A_2,A_3,A_4)$ of commuting bounded Hilbert space operators and a hyperinvariant subspace $\frk{X}_1$ for $A$ such that the Taylor joint spectrum of $A$ restricted to $\frk{X}_1$ is a not a subset of the Taylor joint spectrum of $A$.
\end{abstract}

Let $T=(T_1,\ldots,T_d)$ be a $d$-tuple of commuting bounded operators on a Hilbert space $\frk{H}$, and denote by $\TSpec(T)$ the Taylor joint spectrum \cite{Taylor1970Spec} of $T$.
A subspace $\frk{M}$ of $\frk{H}$ is said to be \textit{hyperinvariant} for $T$ if $B\frk{M}\subset\frk{M}$ for every bounded operator $B$ commuting with $T_1,\ldots,T_d$.
The purpose of this note is to demonstrate that there exist a quadruple $A=(A_1,A_2,A_3,A_4)$ of commuting bounded Hilbert operators and a subspace $\frk{X}_1$ that is hyperinvariant for $A$ such that $\TSpec(A|_{\frk{X}_1})$ is not a subset of $\TSpec(A)$.
The existence of such tuples seems to be folk theorem, but we are motivated by other work to record it.
That $A$ should exist may seem surprising when compared with the single-operator case.
In particular, if $A_0$ is a bounded operator on a Hilbert space $\frk{H}$ and $\frk{M}$ is an $A_0$-hyperinvariant subspace of $\frk{H}$, then $\sigma(A_0|_{\frk{M}})\subset \sigma(A_0)$.

Before detailing the construction of $A$, we would like to thank J. Eschmeier for suggesting that we look at a particular example of J. L. Taylor \cite[Sec. 4]{Taylor1970Spec}.
We would also like to thank R. Clou\^{a}tre for his help in the development of this note.

In what follows, let $\Delta\subset \CC^2$ be a closed polydisc with non-empty interior that contains $0$, and let $U$ be a bounded open neighborhood of $\Delta$.
Set $V=U\bksl \Delta$, which we view alternately as a subset of $\RR^4$.
We denote by $\zeta_1,\zeta_2$ the canonical complex coordinates of $V$ and by $x_1,x_2,x_3,x_4$ the canonical real coordinates, with
\[ \zeta_j=x_{2j-1}+ix_{2j}, \quad j=1,2. \]

Given a positive integer $k$, denote by $C^k_0(V)$ the vector space of compactly supported functions on $V$ having continuous derivatives up to and including order $k$.
For $\phi\in C^1_0(V)$, we set
\[ D_{x_j}\phi =\frac{\pd \phi}{\pd x_j}, \quad j=1,2,3,4. \]
Let $L^2(V)$ denote the $L^2$-space on $V$ with Lebesgue measure $\tau$.
We adopt the convention of identifying an element $f\in L^2(V)$ with a function in the case where $f$ has a smooth representative.

Given $f\in L^2(V)$ and $j\in\{1,2,3,4\}$, recall that $g\in L^2(V)$ is the \textit{weak $x_j$-derivative} of $f$ if
\[ \int_V g\phi d\tau = -\int_V f D_{x_j}\phi d\tau \]
for every $\phi\in C_0^{1}(V)$.
Because $C_0^1(V)$ is dense in $L^2(V)$, the weak $x_j$-derivative of $f$, when it exists, is unique; we denote it by $D_{x_j}f$.
For $\phi\in C^1_0(V)$, the $x_j$-derivative of $\phi$ and the weak $x_j$-derivative of $\phi$ determine the same element of $L^2(V)$, and we denote them both by $D_{x_j}\phi$.
It should be noted that, from our definition, the existence of a weak $x_j$-derivative for $f\in L^2(V)$ implies that $D_{x_j}f\in L^2(V)$.

Denote by $W^{1,2}(V)$ the linear space consisting of all those $f\in L^2(V)$ for which $D_{x_1}f,\ldots,D_{x_4}f$ exist.
Equipped with the norm given by
\[ \|f\|_{W^{1,2}}=\left(\|f\|^2_{L^2(V)}+\sum_{j=1}^4 \|D_{x_j}f\|^2_{L^2(V)}\right)^{1/2}, \]
$W^{1,2}(V)$ is a Hilbert space.
Given $f\in W^{1,2}(V)$, we set
\[ D_{\zeta_j}f=\frac{1}{2}(D_{x_{2j-1}}f-i D_{x_{2j}}f), \quad D_{\cc{\zeta}_j}f=\frac{1}{2}(D_{x_{2j-1}}f+i D_{x_{2j}}f) \]
for $j=1,2$.

A element $f\in W^{1,2}(V)$ is a \textit{weak solution} of the Laplace equation on $V$ if
\[ \sum_{j=1}^4 \int_V (D_{x_j}f)(D_{x_j}\phi)d\tau =0  \]
for every $\phi\in C_0^1(V)$.
Given $f\in W^{1,2}(V)$, 
\begin{equation}\label{Eq1} \sum_{j=1}^4\int_V (D_{x_j}f)(D_{x_j}\phi)d\tau = 4\sum_{j=1}^2\int_V (D_{\cc{\zeta}_j}f)(D_{\zeta_j}\phi)d\tau \end{equation}
when $\phi\in C^2_0(V)$.
Thus, by a simple density argument, \eqref{Eq1} also holds when $\phi\in C_0^1(V)$.
In particular, if $D_{\cc{\zeta}_1}f=D_{\cc{\zeta}_2}f=0$ for some $f\in W^{1,2}(V)$, then $f$ is a weak solution of Laplace's equation on $V$.
It then follows from \cite[Cor. 8.11]{GilTrud_EllipticPDE} that $f$ is infinitely differentiable.
Because $D_{\cc{\zeta}_1}f=D_{\cc{\zeta}_2}f=0$, we see that $f$ is analytic on $V$.

Set $\frk{X}=W^{1,2}(V)\oplus L^2(V)$, and define $A_1,A_2,A_3,$ and $A_4$ in $\mc{B}(\frk{X})$ as follows.
Given $(f,g)\in\frk{X}$, set
\[	A_1(f,g) = (\zeta_1f,\zeta_1 g), \quad	  A_2(f,g) = (\zeta_2f,\zeta_2 g), \]
\[	A_3(f,g) = (0, D_{\cc{\zeta}_1}f), \quad	A_4(f,g) = (0, D_{\cc{\zeta}_2}f), \]
and set $A=(A_1,A_2,A_3,A_4)$.
It is easily seen that $A_1,A_2,A_3,$ and $A_4$ commute.
Set $\frk{X}_1=\ker(A_3)\cap\ker(A_4)$.
It should be noted that if $(f,g)\in\frk{X}_1$, then $D_{\cc{\zeta}_1f}=D_{\cc{\zeta}_2}f=0$ and thus $f$ is analytic.
We also note that $\frk{X}_1$ is a hyperinvariant subspace for $A$.

\vspace{10pt}
\noindent\textbf{Theorem.} \textit{$\TSpec(A|_{\frk{X}_1})$ is not a subset of $\TSpec(A)$.}

\begin{proof}
	Given a bounded $C^\infty$-function $\psi$ on $V$, we note that $M_\psi(f,g)=(\psi f,\psi g)$ is a bounded operator on $\frk{X}$.
	Set $\psi_j=\cc{\zeta}_j/(|\zeta_1|^2+|\zeta_2|^2)$ for $j=1,2$, and note that $\psi_1,\psi_2\in C^\infty(V)$ are bounded and satisfy the equation $1=\zeta_1\psi_1+\zeta_2\psi_2$ on $V$.
	Thus $I=M_{\psi_1}A_1+M_{\psi_2}A_2$, and it then follows from \cite[Lemma 1.1]{Taylor1970Spec} and the corollary to \cite[Lemma 1.3]{Taylor1970Spec} that $0\nin\TSpec(A)$.
	Note that, if $R=[A_1|_{\frk{X}_1} \: \ldots \: A_4|_{\frk{X}_1}]$ from $\frk{X}_1^{\oplus 4}$ into $\frk{X}_1$ is not surjective, then it follows from \cite[Def. 1.1]{Taylor1970Spec} that $0\in \TSpec(A|_{\frk{X}_1})$.

	Suppose that $R$ is surjective.
	Then there would exist $(f_1,g_1),\ldots,(f_4,g_4)\in\frk{X}_1$ such that
	\[ \sum_{j=1}^4 A_j(f_j,g_j)=(1,0). \]
	This implies, in particular, that
	\[ \zeta_1f_1+\zeta_2f_2=1 \]
	where $f_1,f_2$ are analytic.
	By Hartogs's Theorem \cite[Theorem I.C.5]{GnR}, $f_1,f_2$ would have unique analytic extensions $\widetilde{f}_1,\widetilde{f}_2$ to $U$ such that
	\[ 1=z_1\widetilde{f}_1(z_1,z_2)+z_2\widetilde{f}_2(z_1,z_2), \quad (z_1,z_2)\in U. \]
	However, $(0,0)$ is in $U$.
	Thus $R$ cannot be surjective and thus $0\in\TSpec(A|_{\frk{X}_1})$.
\end{proof}

\bibliographystyle{plain}
\bibliography{biblio_main.bib}

\end{document}